\documentclass[reqno,11pt]{amsart} 
\usepackage{amsmath}
\usepackage{amssymb}
\usepackage{amsthm}
\usepackage{verbatim}
\usepackage{xcolor}
\usepackage{graphics}

\newcommand\NoBlackBoxes{\global\overfullrule0pt}

\NoBlackBoxes
\theoremstyle{plain}

\begin{document}

\title{WEIGHTED CKP INEQUALITIES INVOLVING \\
R\'ENYI DIVERGENCE POWERS
}

\author{Sergey G. Bobkov$^{1}$}
\thanks{1) 
School of Mathematics, University of Minnesota, Minneapolis, MN, USA,
bobkov@math.umn.edu. 
}

\author{Devraj Duggal$^{2}$}
\thanks{2) School of Mathematics, University of Minnesota, Minneapolis, MN, USA,
dugga079@umn.edu.
}


\begin{abstract}
Pinsker-type inequalities are considered for the weighted total variation distance between 
probability measures in terms of the R\'enyi divergence powers. They are applied in derivation 
of transport-entropy inequalities under moment-type conditions.
\end{abstract}

\maketitle
\markboth{Sergey G. Bobkov and Devraj Duggal}{Weighted CKP inequalities}

\def\theequation{\thesection.\arabic{equation}}
\def\E{{\mathbb E}}
\def\R{{\mathbb R}}
\def\C{{\mathbb C}}
\def\P{{\mathbb P}}
\def\Z{{\mathbb Z}}
\def\L{{\mathbb L}}
\def\T{{\mathbb T}}

\def\G{\Gamma}

\def\Ent{{\rm Ent}}
\def\var{{\rm Var}}
\def\Var{{\rm Var}}

\def\H{{\rm H}}
\def\Im{{\rm Im}}
\def\Tr{{\rm Tr}}
\def\s{{\mathfrak s}}

\def\k{{\kappa}}
\def\M{{\cal M}}
\def\Var{{\rm Var}}
\def\Ent{{\rm Ent}}
\def\O{{\rm Osc}_\mu}

\def\ep{\varepsilon}
\def\phi{\varphi}
\def\vp{\varphi}
\def\F{{\cal F}}

\def\be{\begin{equation}}
\def\en{\end{equation}}
\def\bee{\begin{eqnarray*}}
\def\ene{\end{eqnarray*}}

\thispagestyle{empty}

Contents:

\vskip2mm
1. Introduction

2. R\'enyi divergence powers

3. Transport-entropy  inequalities

4. Linearization of R\'enyi divergence power 

5. Uniquence of extremizer and its description

6. H\"older-type inequalities for densities

7. Proof of Theorems 4.1--4.2

8. Necessary and sufficient conditions. Moment bounds

9. Proof of Theorem 2.3 (the case $\alpha \geq 2$)

10. Proof of Theorem 2.3 (the case $1 \leq \alpha \leq 2$)

11. Pearson-Vajda distances. Proof of Theorems 2.1--2.2

\vskip5mm
\section{{\bf Introduction}}
\setcounter{equation}{0}

\vskip2mm
\noindent
Let $(E,\mathcal E,\mu)$ be a probability space, and $\nu$ be a probability measure on $E$
which is absolutely continuous with respect to $\mu$, for short, $\nu <\!\!< \mu$,
with density $f = d\nu/d\mu$.
The Pinsker's, also called Csisz\'ar-Kullback-Pinsker's inequality, connects the total variation norm
$v = \|\nu - \mu\|_{\rm TV}$ between these measures
with the informational divergence (Kullback-Leibler distance)
$$
D = D(\nu||\mu) = \int f\,\log f\,d\mu
$$
by means of the inequality
\be
\|\nu - \mu\|_{\rm TV} \leq \sqrt{2 D}.
\en
Often, the informational divergence is much easier treated in comparison with total variation,
and then the inequality (1.1) is of particular use (for example, when both $\mu$ and $\nu$ 
are Gaussian measures on the Euclidean space).

The square root function appearing on the right-hand side of (1.1) is not optimal.
Various improvements in the form $D \geq \Psi(v)$, which we do not consider here,
have been obtained by many authors, see Kullback \cite{K1}, \cite{K2}, Vajda \cite{V},
Toussant \cite{Tou}, Topsoe \cite{Top}, Fedotov, Harremo\"es and Topsoe \cite{F-H-T}.
Another way of sharpening, important for applications,
is based on the replacement of $v$ with the weighted total variation distance
$$
\|w(\nu - \mu)\|_{\rm TV} = \int w\,d|\nu-\mu| = \int w\,|f-1|\,d\mu,
$$
where $w \geq 0$ is a given measurable function on $E$ (weight function).
However, to get an anolog of (1.1) in this more general setting, one needs to require that 
$w$ has certain integrability properties under the measure $\mu$. Interesting results in this direction
were obtained by Bolley and Villani \cite{B-V} who derived the bounds
\be
\|w (\nu - \mu)\|_{\rm TV} \leq \Big(1 + \log \int e^{w^2} d\mu\Big)^{1/2}
\sqrt{2 D}
\en
and, as another variant,
\be
\|w (\nu - \mu)\|_{\rm TV} \leq \Big(\frac{3}{2} + \log \int e^{2w} d\mu\Big)
\Big(\sqrt{D} + \frac{1}{2}\, D\Big).
\en
Choosing $w=1$, (1.2) yields (1.1) with an additional factor of $\sqrt{2}$. The point of (1.3)
is that a weaker integrability assumption is required from $w$ in comparison with (1.2), 
although this is achieved at the expense of an additional $D$-term in the last brackets.

Bolley and Villani call (1.2)-(1.3) weighted CKP inequalities. They applied these results
to derive transport-entropy bounds in metric measure spaces where the metrics have proper
exponential moments under $\mu$. Later on,
the inequality (1.2) was also used in \cite{B-C-G1} in the study of the central limit
theorem with respect to the relative Fisher information.

\vskip7mm
\section{{\bf R\'enyi Divergence Powers}}
\setcounter{equation}{0}

\vskip2mm
\noindent
The purpose of this paper is to explore weighted CKP inequalities under weaker moment
conditions posed on the weight function $w$. This naturally requires to involve stronger
informational distances in place of $D$ such as the R\'enyi divergence powers, which are also 
known under the name ``Tsallis distances".

Keeping the same setting as before, the R\'enyi divergence of order $\alpha>1$
between $\nu$ and $\mu$ is defined by
$$
D_\alpha = D_\alpha(\nu||\mu) = \frac{1}{\alpha - 1}\,\log\int f^\alpha\,d\mu,
$$
provided that $\nu <\!\!< \mu$ has density $f = d\nu/d\mu$, and $D_\alpha(\nu||\mu) = \infty$
otherwise. This quantity is related to the R\'enyi's entropy like the Kullback-Leibler divergence  $D$
is related to Shannon's entropy. Correspondingly, by analogy with entropy power,
the R\'enyi divergence power or the Tsallis distance is given by
\be
T_\alpha = T_\alpha(\nu||\mu) = \frac{1}{\alpha - 1}\,\left[\int f^\alpha\,d\mu - 1\right].
\en
These quantites are non-decreasing in $\alpha$ with limit $T_1 = D_1 = D$ (once $T_\alpha$ 
is finite for some $\alpha>1$) and are connected by 
$$
T_\alpha = \frac{1}{\alpha - 1}\,\big[e^{(\alpha - 1) D_\alpha} - 1\big] \geq D_\alpha. 
$$
Hence, they are of the same order, when these distances are small, although their role
in various relations might be different. We refer an interested reader to \cite{E-H1}, \cite{E-H2},
\cite{B-C-G2} for an account of basic properties of these informational functionals.

In the sequel, we denote by 
$$
\|u\|_p = \Big(\int |u|^p\,d\mu\Big)^{1/p}
$$ 
the $L^p(\mu)$-norm of a measurable function $u$ on $E$, and by $\beta = \frac{\alpha}{\alpha-1}$ 
the conjugate power.

We aim to derive several CKP bounds on the weighted total variation distance in terms 
of the $L^\beta$-norms of the weighed function $w$ and the R\'enyi divergence power as defined in (2.1). 
Recall that the probability measure $\mu$ on $E$ is fixed, while $\nu$ may be arbitrary.

\vskip5mm
{\bf Theorem 2.1.} {\sl If $1 < \alpha \leq 2$, then
\be
\|w (\nu - \mu)\|_{\rm TV} \leq \frac{16}{3}\,\|w\|_\beta\,
\max\big\{\sqrt{T_\alpha}, T_\alpha^{1/\alpha}\big\}.
\en
If $\alpha \geq 2$, then
\be
\|w (\nu - \mu)\|_{\rm TV}  \leq 3^\alpha\,\|w\|_\beta\,T_\alpha^{1/\alpha}.
\en
}

\vskip5mm
As an alternative variant, we also have:

\vskip5mm
{\bf Theorem 2.2.} {\sl Putting $\beta^* = \max(\beta,2)$, 
\be
\|w (\nu - \mu)\|_{\rm TV} \leq C_\beta\,\|w\|_{2\beta^* - 2}\,\sqrt{T_\alpha}.
\en
Here one may take $C_\beta = 2$ for $\alpha \geq 2$ and $C_\beta =  4\,\beta^\beta$
for $\alpha<2$.
}

\vskip5mm
Note that in the limit case $\alpha \rightarrow 1$, the inequalitiy (2.2) with $w=1$ 
returns us to the Pinsker inequality (1.1) with an additional numerical factor.

The proof Theorem 2.1 is rather simple and involves an application of Pearson-Vajda 
informational distances which are closely related to $T_\alpha$ (cf. Section 11).
As for Theorem 2.2, the argument is based on the linearization of the R\'enyi divergence power developed
in \cite{B-D}. Here we will recall and refine this argument, and provide some missing technical details.
Another intermediate step in derivation of (2.4) deals with relations of the form
\be
\Big|\int u\,d\nu\Big| \leq K\sqrt{T_\alpha}, 
\en
which are required to hold for a fixed function $u \in L^1(\mu)$ with $\mu$-mean zero 
in the class of all probability measures $\nu <\!\! < \mu$. An interesting question of independent
interest is how one can describe or estimate the best constant $K$ in (2.5). The next assertion
which was essentially obtained, but not properly emphasized in \cite{B-D} gives an answer in terms
of the quantities $K_p = K_p(u)$ defined by
\be
K_p^{2p-2} \, = \, 
\sup_{r > 0} \, \Big[r^{p - 2} \int_{|u| \geq r} |u|^p\,d\mu\Big], \quad p \geq 2.
\en
In particular,  $K_2 = \|u\|_2$ is the usual $L^2$-norm.

\vskip5mm
{\bf Theorem 2.3.} {\sl $a)$ \ Let $1 < \alpha \leq 2$. The best value of $K$ in $(2.5)$
satisfies
\be
c_\beta K_{\beta} \leq K \leq C_\beta K_{\beta},
\en
holding with $c_\beta = \frac{1}{4}\,\beta^{-\beta}$ and $C_\beta = 2\beta^\beta$.

$b)$ For $\alpha \geq 2$, the best value of $K$ satisfies 
$$
\sqrt{\frac{2}{\alpha}}\,\|u\|_2 \leq K \leq \|u\|_2.
$$
}

\vskip2mm
For completeness, we will inlcude the proof of this theorem (polishing 
some steps in the original argument).
Note that one may relate $K_p$ to the Lebesgue norms, assuming
that $\|u\|_{2p-2}$ is finite. In this case, by Markov's inequality, for all $r>0$, 
$$
 \int_{|u| \geq r} |u|^p\,d\mu \leq \frac{1}{r^{p-2}}\,\|u\|_{2p-2}^{2p-2}
$$
and hence (2.6) yields
\be
K_p(u) \leq \|u\|_{2p-2}.
\en

\vskip5mm
{\bf Corollary 2.4.} {\sl The best value of $K$ in the inequality $(2.5)$ satisfies
\be
K \leq C_\beta\,\|u\|_{2\beta^* - 2}, \quad \beta^* = \max(\beta,2).
\en
Here one may take $C_\beta = 1$ for $\alpha \geq 2$ and $C_\beta =  2\beta^\beta$
for $\alpha<2$.
}

\vskip7mm
\section{{\bf Transport-Entropy Inequalities}}
\setcounter{equation}{0}

\vskip2mm
\noindent
Theorems 2.1-2.2 can be applied in the derivation of transport-entropy inequalities.
Such inequalities provide upper bounds for the Kantorovich distances between
probability distributions in terms of informaton-theoretic distances. 

Let $(E,\rho)$ be a separable metric space endowed with the $\sigma$-algebra $\mathcal E$ 
of Borel subsets in $E$. The Kantorovich transport distance of order $p \geq 1$ 
between Borel probability measures $\mu$ and $\nu$ on $E$ is defined by
$$
W_p(\mu,\nu)  =  \Big(\inf \int\!\!\int \rho(x,y)^p\,d\pi(x,y)\Big)^{1/p},
$$
where the infimum is running over all Borel probability measures $\pi$ on the product
space $E \times E$ with marginals $\mu$ and $\nu$.
As is well-known, $W_p$ represents a metric in the space ${\mathcal P}_p(E)$ of all
Borel probability measures $\mu$ on $E$ such that
\be
M_p = \Big(\int \rho(x,x_0)^p\,d\mu(x)\Big)^{1/p} < \infty
\en
for some point $x_0 \in E$ (equivalently, for all $x_0$).  See e.g. \cite{V}.

If $p=1$, the Kantorovich duality theorem asserts that
$$
W_1(\mu,\nu) = \sup \Big|\int u\,d\nu - \int u\,d\mu \Big|,
$$
where the infimum is taken over all functions $u$ on $E$ with Lipschitz
semi-norm $\|u\|_{\rm Lip} \leq 1$, cf. \cite{D}. Using this description together with
Theorem 2.3, it was shown in \cite{B-D} that the transport-entropy inequality
\be
W_1(\mu,\nu) \leq K \sqrt{T_\alpha}, \quad \nu \in {\mathcal P}_1(E),
\en
holds true for a given measure $\mu$ in ${\mathcal P}_1(E)$ with some finite 
constant $K$, if and only if
$$
\sup_{r > 0} \, \Big[r^{\beta^* - 2} \int_{\rho(x,x_0) \geq r} 
\rho(x,x_0)^{\beta^*}\,d\mu\Big] < \infty, \quad \beta^* = \max(\beta,2).
$$
In particular, by (2.8), this condition is fulfilled as long as (3.1) holds
with $p = 2\beta^* - 2$.

In the limit case as $\alpha \rightarrow 1$, the relation (3.2) becomes
\be
W_1(\mu,\nu) \leq K \sqrt{D}, \quad \nu \in {\mathcal P}_1(E),
\en
which was previously studied in \cite{B-G}. However, the moment condition (3.1)
is unsufficient for the finiteness of the constant $K$, since (3.3) is equivalent
$$
\int e^{c \rho(x,x_0)^2}\,d\mu(x) < \infty
$$
with some $c>0$. This is also equivalent to the property that all Lipschitz functions on $E$
are subgaussian under $\mu$.

Hence the subgaussianity property is also necessary when $W_1$ is replaced in (3.3) with a stronger 
metric $W_p$ for $p>1$. In fact, this case turns out to be rather different and more difficult. 
As an important example, when $\mu = \gamma_n$ is the standard Gaussian 
measure on the Euclidean space $E = \R^n$, there is a remarkable inequality due 
to Talagrand \cite{T} (cf, also \cite{B-G}):
$$
W_2(\gamma_n,\nu) \leq \sqrt{2 D(\nu||\gamma_n)}.
$$
A similar relation such as
\be
W_2(\mu,\nu) \leq K\sqrt{D(\nu||\mu)}
\en
remains to hold for probability measures $\mu$ on $\R^n$ having a log-concave density
with respect to $\gamma_n$, even if after some rescaling (\cite{O-V}). Let us refer to 
\cite{G}, \cite{G-L} for more results on the subject and only note that
a full characterization of the class of measures $\mu$ on $\R^n$ for which (3.4) holds true
with some constant $K$ is still unknown.

Some generalizations of (3.4) were considered by Bolley and Villani \cite{B-V}.
In the abstract metric space setting,
they employed the following upper bound on $W_p$ in terms of the weighted
total variation distance (\cite{V}, Proposition 7.10):
\be
W_p^p(\mu,\nu) \leq 2^{p-1}\,\|\rho(x,x_0)^p\,(\nu-\mu)\|_{\rm TV}.
\en
Assuming that
\be
\int e^{c \rho(x,x_0)^{2p}}\,d\mu(x) < \infty \quad {\rm or} \quad
\int e^{c \rho(x,x_0)^p}\,d\mu(x) < \infty
\en
for some $c>0$ and applying their general bounds (1.2)-(1.3), 
they respectively derived the transport-entropy inequalities
\bee
W_p(\mu,\nu) 
 & \leq &
K \, D(\nu||\mu)^{\frac{1}{2p}}, \\
W_p(\mu,\nu) 
 & \leq &
K\,\big(D(\nu||\mu)^{\frac{1}{2p}} + D(\nu||\mu)^{\frac{1}{p}}\big)
\ene
with constants $K$ independent of $\nu$.

We can now follow a similar argument so as to replace (3.6) with a  weaker
moment assumption (3.1). This may be achieved at the expense of the replacement of 
the Kullback-Leibler distance with the R\'enyi divergence power.  Applying Theorem 2.1 
with $w(x) = \rho(x,x_0)^p$ and using the notation $T_\alpha = T_\alpha(\nu||\mu)$, we get:

\vskip5mm
{\bf Corollary 3.1.} {\sl If $1 < \alpha \leq 2$, then, for any Borel probability measure $\nu$
on $E$,
$$
W_p(\mu,\nu)  \leq CM_{\beta p}\,
\max\Big\{T_\alpha^{\frac{1}{2p}}, T_\alpha^{\frac{1}{\alpha p}}\Big\}
$$
with some absolute constant $C$. If $\alpha \geq 2$, then
$$
W_p(\mu,\nu) \,  \leq \, 2 \cdot 3^{\alpha/p}\,M_{\beta p}\ T_\alpha^{\frac{1}{\alpha p}}.
$$
}

\vskip3mm
As an alternative variant, one may appeal to Theorem 2.2.

\vskip5mm
{\bf Corollary 3.2.} {\sl Putting $\beta^* = \max(\beta,2)$ for $\alpha > 1$, we have
$$
W_p(\mu,\nu) \,  \leq \,  C_\beta^{1/p}\,M_{(2\beta^* - 2)p}\,T_\alpha^{\frac{1}{2p}}.
$$
One may take $C_\beta = 2$ for $\alpha \geq 2$ and $C_\beta =  16\,\beta^\beta$
for $\alpha<2$.
}

\vskip7mm
\section{{\bf Linearization of R\'enyi Divergence Power }}
\setcounter{equation}{0}

\vskip2mm
\noindent
Let us return to the setting of an abstract probability space $(E,\mathcal E,\mu)$.
Towards the proof of Theorems 2.1-2.2, we recall the basic linearization argument.

Denote by ${\mathcal P}(\mu)$ the collection of all probability densities $f$ on 
$E$ with respect to $\mu$. Being convex on ${\mathcal P}(\mu)$, the entropy functional $D$ 
admits a well-known sup-linear representation
$$
D(\nu||\mu) = \sup \left\{\int fg\,d\mu: \int e^{g}\,d\mu \leq 1\right\}.
$$
In other words,
$$
\int g\,d\nu \leq D(\nu||\mu)
$$
for all $\nu <\!\!< \mu$, if and only if $\int e^{g}\,d\mu \leq 1$. There is a similar description
for the R\'enyi divergence power of an arbitrary order $\alpha>1$. 
In the sequel, we write $t_+ = \max(t,0)$, $t \in \R$.
In a slightly weaker form, the next assertion is proved in \cite{B-D}.

\vskip5mm
{\bf Theorem 4.1.} {\sl The relation
\be
\int g\,d\nu \leq T_\alpha(\nu||\mu)
\en
holds true for any probability measure $\nu <\!\!< \mu$ on $E$ such that the
integral in $(4.1)$ exists, if and only if $g_+ \in L^\beta(\mu)$, and
\be
\int (g-c)_+^\beta\,d\mu \, \leq \, -\beta^\beta\,
\big(c + \beta - 1\big),
\en
where $c$ is a unique solution to the equation
\be
\int (g-c)_+^{\beta - 1}\,d\mu \, = \, \beta^{\beta - 1}.
\en
}

\vskip3mm
As an equivalent description, Theorem 4.1 admits the following analog.

\vskip5mm
{\bf Theorem 4.2.} {\sl The relation $(4.1)$ holds true for any probability measure 
$\nu <\!\!< \mu$ on $E$  such that the integral in $(4.1)$ exists, if and only if 
$g_+ \in L^\beta(\mu)$, and the condition $(4.2)$ is fulfilled for at least one constant $c$.
}

\vskip5mm
In \cite{B-D} Theorems 4.1-4.2 are stated with the additional assumption that $g_+ \in L^\beta(\mu)$.
In the current formulation, this property appears as a necessary condition for (4.1).

Note that the integral in $(4.1)$ exists, althought it might be equal to $-\infty$, if
$g_+ \in L^\beta(\mu)$ and the distance $T_\alpha(\nu||\mu)$ is finite. Indeed,
the latter means that the measure $\nu$ is absolutely continuous with respect to $\mu$ 
and its density $f = d\nu/d\mu$ belongs to $L^\alpha(\mu)$. In this case, by H\"older's inequality,
$$
\int g_+\,d\nu  = \int fg_+\,d\mu \leq \|f\|_\alpha \|g_+\|_\beta < \infty.
$$
Hence, the relation (4.1) is well-defined. Of course, when $T_\alpha(\nu||\mu) = \infty$,
there is no reason to worry about the existence of the integral in (4.1).

On ${\mathcal P}(\mu)$ introduce the concave functional
\be
R f \, = \, \int fg\,d\mu - T_\alpha(\nu||\mu) \, = \, 
\int fg\,d\mu - \frac{1}{\alpha - 1}\,\left[\int f^\alpha\,d\mu - 1\right]
\en
with convention that $R f = -\infty$ in the case $\int f^\alpha\,d\mu = \infty$.
Thus, Theorem 4.1 gives the description of probability measures $\mu$ on $E$ such that
$Rf \leq 0$ for all densities $f$. 

One natural approach to obtain such a description is to identify the maximizer $f$ for 
the functional $R$. Let us recall the proof of Theorems 4.1-4.2 and give some more technical
details, which are missing in \cite{B-D}. This will be the subject of the current and next two sections.
The argument is based on several preliminary lemmas. 
The main observation on the existence of an extremizer for the functional $R$ is the following:

\vskip5mm
{\bf Lemma 4.3.} {\sl If $g_+ \in L^\beta(\mu)$,
the functional $R f$ is bounded above on ${\mathcal P}(\mu)$ and attains maximum at some 
function $f \in {\mathcal P}(\mu) \cap L^\alpha(\mu)$. Moreover,
\be
\int f |g|\,d\mu < \infty.
\en
}

\vskip2mm
{\bf Proof.} By H\"older's inequality,
$$
Rf \, \leq \, \|f\|_\alpha \|g_+\|_\beta - 
\frac{1}{\alpha - 1}\,\Big[\|f\|^{\alpha}_\alpha - 1\Big] \, \leq \,
c_0 + c_1 \|g_+\|_\beta^\beta
$$
up to some constants $c_0$ and $c_1$ depending on $\alpha$ and $\|g_+\|_\beta$. Hence
\be
M \, = \, \sup_{f  \in {\mathcal P}(\mu)} R f \, = \,
\sup\big\{Rf: f \in {\mathcal P}(\mu) \cap L^\alpha(\mu)\big\} < \infty.
\en

For the normlized indictor functions $f_A = \frac{1}{\mu(A)}\,1_A$ with
$A \in \mathcal E$, $\mu(A) > 0$, we have
$$
R f_A =  \frac{1}{\mu(A)}\,\int_A g\,d\mu - \frac{1}{\alpha - 1}\,\left[ \frac{1}{\mu(A)^{\alpha - 1}} - 1\right].
$$
In particular, choosing $A_c = \{g \geq c\}$ with any parameter $c \in \R$ such that $\mu(A_c) > 0$ leads to
$$
R f_{A_c} \geq  c - \frac{1}{\alpha - 1}\,\left[ \frac{1}{\mu(A_c)^{\alpha - 1}} - 1\right] > -\infty.
$$
This shows that $M > -\infty$. Thus, the supremum in (4.6) is finite. 

Moreover, one may restrict this supremum to all $f$ such that that $\|f\|_{\alpha} \leq C$ 
with some large $C$. Indeed, if $\|f\|_{\alpha} > C$, the expression
$$
\|f\|_\alpha \|g_+\|_\beta - 
\frac{1}{\alpha - 1}\,\Big[\|f\|^{\alpha}_\alpha - 1\Big]
$$
tends to $-\infty$ for $C \rightarrow \infty$.
Therefore, $R$ is bounded above on $\mathcal P$ by the finite constant 
$$
M \, = \, 
\sup\big\{Rf: f \in {\mathcal P}(\mu) \cap L^\alpha(\mu), \ \|f\|_{\alpha} \leq C\big\}
$$
for a sufficiently large value of $C$.

Take a sequence $f_n \in {\mathcal P}(\mu) \cap L^\alpha(\mu)$ with 
$\|f_n\|_{\alpha} \leq C$ such that $R f_n \ \rightarrow M$ as $n \rightarrow \infty$.
The unit ball of $L^\alpha$ is weakly compact, so there is a subsequence
$f_{n'}$ weakly convergent to some $f$ with $\|f\|_{\alpha} \leq C$. For simplicity, let
this subsequence be the whole sequence, so that
\be
\int f_n h\,d\mu \rightarrow \int f h\, d\mu \quad {\rm as} \ n \rightarrow \infty
\en
for all $h \in L^\beta(\mu)$. Necessarily $f \in {\mathcal P}(\mu)$. Indeed, being applied to indicator
functions $h = 1_A$, $A \in \mathcal E$, (4.7) implies that $\int_A f\,d\mu \geq 0$.
Since $A$ is arbitrary measurable, it follows that $f \geq 0$ $\mu$-a.e. Being applied
with $h = 1$, we also obtain that $\int f\,d\mu = 1$, hence $f$ is a probability density on $E$
with respect to $\mu$. It is also well-known that the weak convergence implies that
\be
\|f\|_{\alpha} \, \leq \, \liminf_{n \rightarrow \infty}\, \|f_n\|_\alpha.
\en

Now, since $g_+ \in L^\beta(\mu)$, one may apply (4.7) to $h = g_c = \max(g,c)$ with 
parameter $c \in \R$ (which will be sent to $-\infty$). Using $g \leq g_c$, we have
$f_n g \leq f_n g_c$ $\mu$-a.e., and (4.7) yields
$$
\limsup_{n \rightarrow \infty} \int f_n g\,d\mu \, \leq \,
\limsup_{n \rightarrow \infty} \int f_n g_c\,d\mu \, = \, \int f g_c\,d\mu.
$$
Applying also (4.8), we thus get
\bee
\limsup_{n \rightarrow \infty}\, Rf_n
 & = &
\limsup_{n \rightarrow \infty}\, \Big(\int f_n g\,d\mu - 
\frac{1}{\alpha - 1}\,\big[\|f_n\|_\alpha^\alpha - 1\big]\Big) \\
 & \leq &
\limsup_{n \rightarrow \infty}\,\int f_n g\,d\mu +
\limsup_{n \rightarrow \infty}\,\Big(-  \frac{1}{\alpha - 1}\,\big[\|f_n\|_\alpha^\alpha - 1\big]\Big) \\
 & \leq  &
\int f g_c\,d\mu -  \frac{1}{\alpha - 1}\, 
\liminf_{n \rightarrow \infty}\big[\|f_n\|_\alpha^\alpha - 1\big]\\
 & \leq &
\int fg_c\,d\mu - \frac{1}{\alpha - 1}\,\big[\|f\|_\alpha^\alpha - 1\big].
\ene
As a result,
$$
\limsup_{n \rightarrow \infty}\, Rf_n \leq \int fg_c\,d\mu - 
\frac{1}{\alpha - 1}\,\big[\|f\|_\alpha^\alpha - 1\big].
$$
Here, the left-hand side does not depend on $c$, while 
$\int fg_c\,d\mu \downarrow \int fg\,d\mu$ as $c \rightarrow -\infty$, by the monotone
convergence theorem. Hence, in the limit we get
$$
\limsup_{n \rightarrow \infty}\, Rf_n \leq \int fg\,d\mu - 
\frac{1}{\alpha - 1}\,\big[\|f\|_\alpha^\alpha - 1\big] = Rf.
$$
It follows that $R f = M$, which proves the first claim.

As was already noticed, $\int f g_+\,d\mu < \infty$. Since $M$ is finite and $\|f\|_{\alpha} \leq C$,
necessarily $\int f g\,d\mu > -\infty$, so that $\int f g_-\,d\mu < \infty$ for the function $g_- = \max(-g,0)$.
This proves the second claim (4.5).
\qed

\vskip7mm
\section{{\bf Uniquence of Extremizer and its Description}}
\setcounter{equation}{0}

\vskip2mm
\noindent
The next preliminary step about the functional $R$ defined in (4.4) is the following
lemma from \cite{B-D}, whose proof requires some missing technical step.
Speaking about the uniqueness, we identify the functions in $\mathcal P(\mu)$ that coincide
$\mu$-almost everywhere.

\vskip5mm
{\bf Lemma 5.1.} {\sl If $g_+\! \in L^\beta(\mu)$, the maximizer for the functional 
$R$ is unique and has the form
\be
f \, = \, \beta^{1-\beta}\, (g-c)_+^{\beta - 1} \quad \mu-{\sl a.e.}
\en
for some (unique) constant $c$.
}

\vskip5mm
{\bf Proof.} Let $f \in {\mathcal P}(\mu) \cap L^\alpha(\mu)$ be a maximizer for $R$ as in Lemma 4.3. 
For $\delta > 0$, put $A_\delta = \{x \in E: f(x) > \delta\}$. Since $f \geq 0$ and 
$\int f\,d\mu = 1$, we have $\mu(A_\delta) > 0$ for all $\delta > 0$ small enough. This will be assumed.

Consider the functions of the form
\be
f_\ep = f + \ep u, \quad \ep \in \R,
\en
where $u$ is a bounded measurable function on $E$ vanishing outside $A_\delta$
and such that $\int u\,d\mu = 0$ and $\int |u g|\,d\mu < \infty$. Then, $f_\ep$ will belong 
to ${\mathcal P}(\mu) \cap L^\alpha(\mu)$ for all sufficiently small $\ep$ and hence
$Rf_\ep \leq Rf$.

Using Taylor's expansion, one can show that, as $\ep \rightarrow 0$,
\be
\int f_\ep^\alpha\,d\mu = \int f^\alpha\,d\mu + 
\alpha \ep \int f^{\alpha - 1} u\,d\mu + o(\ep).
\en
Hence, by (4.4),
\bee
Rf_\ep 
 & = &
\int fg\,d\mu - \frac{1}{\alpha - 1}\,\left[\int f^\alpha\,d\mu - 1\right] \\
 & & +  \ \ep\,
\Big(\int u g\,d\mu - \frac{\alpha}{\alpha - 1}\, \int f^{\alpha - 1} u\,d\mu\Big) + o(\ep).
\ene
Therefore,
$$
Rf_\ep - Rf = \ep \int \big(g - \beta f^{\alpha - 1}\big)\,u\,d\mu + o(\ep).
$$
Since $\ep$ may be both positive and negative (although small), from $Rf_\ep \leq Rf$
it follows that
$$
\int \big(g - \beta f^{\alpha - 1}\big)\,u\,d\mu = 0
$$
for all admissible functions $u$. This holds, in particular,  for all bounded measurable 
functions $u$ on $E$ vanishing outside the set 
$$
A_{\delta.N} = \{x \in E: f(x) > \delta, \ g(x) \geq -N\}, \quad N = 1,2,\dots, 
$$
But this is only possible when
$g - \beta f^{\alpha - 1} = c$ on $A_{\delta,N}$ $\mu$-a.e. for some constant $c$. 
Since $\delta > 0$ may be arbitrary (although small) and $N$ may also be arbitrary, 
this constant cannot depend on $\delta$ and $N$. As a result, 
\be
g - \beta f^{\alpha - 1} = c \quad \mu-{\rm a.e.} \ {\rm on \ the \ set} \ A_0 = \{x \in E: f(x) > 0\}. 
\en

\vskip2mm
{\bf Case 1:} $\mu(A_0) = 1$. Since
$\frac{1}{\alpha - 1} = \beta - 1$, the function $f$ has the stated form (5.1).

\vskip2mm
{\bf Case 2:} $0 < \mu(A_0) < 1$. This case was not considered in \cite{B-D}, so we need
to add missing arguments. Let us return to the definition (5.2), assuming that
$u$ is bounded, has $\mu$-mean zero, and $\int |u g|\,d\mu < \infty$. 
We also need to require that $u \geq 0$ on the set
$$
\bar A_0 = \{x \in E: f(x) = 0\},
$$ 
with $\ep > 0$ being small enough, so that $f_\ep \geq 0$ $\mu$-a.e. and $\int f\,d\mu = 1$. 

If additionally $u$ is vanishing on the set $A_0 \setminus A_\delta$ for some $\delta>0$, we have the
Taylor's expansion
$$
\int_{A_0} f_\ep^\alpha\,d\mu = \int_{A_0} f^\alpha\,d\mu + 
\alpha \ep \int_{A_0} f^{\alpha - 1} u\,d\mu + o(\ep),
$$
and since $f = \ep u$ on $\bar A_0$,
$$
\int_{\bar A_0} f_\ep^\alpha\,d\mu = \ep^\alpha \int_{\bar A_0} u^\alpha\,d\mu.
$$
Adding the two representations, we arrive similarly to (5.3) at the representation
$$
\int f_\ep^\alpha\,d\mu = \int f^\alpha\,d\mu + 
\alpha \ep \int_{A_0} f^{\alpha - 1} u\,d\mu + o(\ep).
$$
Thus,
$$
Rf_\ep \, = \, Rf +  \ep\,
\Big(\int u g\,d\mu - \frac{\alpha}{\alpha - 1}\, \int_{A_0} f^{\alpha - 1} u\,d\mu\Big) 
+ o(\ep),
$$
that is,
$$
Rf_\ep - Rf = \ep \int_{A_0} \big(g - \beta f^{\alpha - 1}\big)\,u\,d\mu + 
\ep \int_{\bar A_0} u g\,d\mu  + o(\ep).
$$
Using the relation (5.4) for the first integrand, the above is simplified as
$$
Rf_\ep - Rf = \ep \int_{A_0}c\,u\,d\mu + 
\ep \int_{\bar A_0} u g\,d\mu  + o(\ep),
$$
or equivalently (since $u$ has $\mu$-mean zero)
$$
Rf_\ep - Rf = \ep \int_{\bar A_0} u (g-c)\,d\mu  + o(\ep).
$$
Diving this by $\ep$ and letting $\ep \rightarrow 0$, the property $R_\ep f \leq Rf$ leads to
\be
\int_{\bar A_0} u (g-c)\,d\mu \leq 0.
\en

On this step, the assumption that $u$ is vanishing on the set $A_0 \setminus A_\delta$ 
may be removed. Hence, (5.5) holds true for any bounded measurable function $u \geq 0$ on
$\bar A_0$ such that $\int_{\bar A_0} |ug|\,d\mu < \infty$. Indeed, having such a function $u$, one
can extend it to the whole space $E$ (that is, one can redifine it on the set $A_0$) in a way so that
the extended function is bounded, has $\mu$-mean zero, and satisfies the integrablity requirement
$\int |ug|\,d\mu < \infty$. As a consequence, (5.5) is equivalent to the property that
$g-c \leq 0$ $\mu$-a.e. on the set $\bar A_0$, and therefore the required equality (5.1), namely
$f \, = \, \beta^{1-\beta}\, (g-c)_+^{\beta - 1}$ $\mu$-a.e., remains to hold on $\bar A_0$ as well.

The two cases yield (5.1) in general. 
Finally, let us see that the constant $c$ is uniquely determined by 
the condition $\int f\,d\mu = 1$. Define
$$
\varphi(c) \, = \, \int (g-c)_+^{\beta - 1}\,d\mu.
$$
This function is continuous, non-increasing, and convex on the real line
with $\varphi(-\infty) = \infty$, $\varphi(\infty) = 0$. Hence, it is (strictly) decreasing for
$c \leq c_0 = {\rm ess\,sup}\ g$, and we have $\varphi(c_0) = 0$. 
In particular, for any $b>0$, the equation $\varphi(c) = b$ has a unique solution $c$.
\qed

\vskip7mm
\section{{\bf H\"older-type Inequalities for Densities}}
\setcounter{equation}{0}

\vskip2mm
\noindent
In order to show that the condition $g_+ \in L^\beta(\mu)$ is necessary for the
relation (4.1),  we prove:

\vskip5mm
{\bf Proposition 6.1.} {\sl Let $g$ be a non-negative measurable function on $E$. The inequality
\be
\int fg\,d\mu \leq K \int f^\alpha\,d\mu
\en
holds with some constant $K = K(g)$ for any probability density $f$ on $E$ with 
respect to $\mu$,
if and only if $g \in L^\beta(\mu)$. In this case, the optimal value of this constant satisfies
\be
\frac{1}{e \alpha}\, \|g\|_\beta \leq K \leq \|g\|_\beta.
\en
}

\vskip3mm
With proper changes, one can drop the assumption that the function $g$ is non-negative in (6.1).
The extreme case where $g$ is non-positive is not interesting; then (6.1) holds with $K=0$.
So, let $g$ be a measurable function on $E$ such that $\mu\{g \geq 0\} > 0$. 

\vskip5mm
{\bf Corollary 6.2.} {\sl The inequality $(6.1)$ holds with some constant $K = K(g)$ 
for any probability density $f$ on $E$ with respect to $\mu$ such that the
integral $\int fg\,d\mu$ exists, if and only if $g_+ \in L^\beta(\mu)$. In this case, the optimal value 
of this constant satisfies
\be
\frac{1}{e\alpha}\, \left(\frac{1}{\mu\{g \geq 0\}}
\int_{g \geq 0} g^\beta\,d\mu\right)^{1/\beta} \leq K \leq \left(\frac{1}{\mu\{g \geq 0\}}
\int_{g \geq 0} g^\beta\,d\mu\right)^{1/\beta}.
\en
}

\vskip5mm
Indeed, (6.1) is reduced to the case where the functions $f$ are supported
on the set $A = \{g \geq 0\}$, and then this inequality may be rewritten as a 
relation (6.1) for the normalized restriction $\mu_A$ of $\mu$ to this set. It remains
to apply Proposition 6.1 to the probability space $(E,\mu_A)$ with the function $g$
restricted to $A$.

As for the proof of Proposition 6.1, we need to involve another functional
\be
Q f = \int fg\,d\mu - \frac{1}{\alpha} \int f^\alpha\,d\mu
\en
on the same space of probability densities $f$, where $g$ is a non-negative
measurable function on $E$. Lemmas 4.3 and 5.1 have the corresponding counteparts.

\vskip5mm
{\bf Lemma 6.3.} {\sl If $g \in L^\beta(\mu)$ is non-negative, the functional $Q f$ is bounded 
above on ${\mathcal P}(\mu) \cap L^\alpha(\mu)$ and attains maximum at some function $f$ 
in this set. Moreover, the maximizer is unique and has the form 
$$
f = (g-c)_+^{\beta - 1} \ \ \mu-{\rm a.e.}
$$
for some (unique) constant $c$.
}

\vskip3mm
By H\"older's inequality,
$$
Qf \, \leq \, \|f\|_\alpha \|g\|_\beta - \frac{1}{\alpha}\,\|f\|^{\alpha}_\alpha \, \leq \,
c_0 + c_1 \|g\|_\beta^\beta
$$
up to some constants $c_0$ and $c_1$ depending on $\alpha$ and $\|g\|_\beta$. Hence
$$
M \, = \, 
\sup\big\{Qf: f \in {\mathcal P}(\mu) \cap L^\alpha(\mu)\big\} < \infty.
$$
The proof that $M$ is attained is identical to the one of Lemma 4.3. The uniqueness part
also follows the same arguments as in the proof of Lemma 5.1.

\vskip5mm
{\bf Proof of Proposition 6.1.} In one direction, if $g \in L^\beta(\mu)$, one may apply
H\"older's inequality to get
$$
\int fg\,d\mu \leq \|g\|_\beta \|f\|_\alpha \leq  \|g\|_\beta \|f\|_\alpha^\alpha,
$$
where we used $1 = \|f\|_1 \leq \|f\|_\alpha \leq \|f\|_\alpha^\alpha$. This gives the upper 
bound on $K$ in (6.2).

For the opposite direction, first assume that $g \in L^\beta(\mu)$ and that (6.1) holds true 
with constant $K = 1/\alpha$, i.e. 
\be
\int fg\,d\mu \leq \frac{1}{\alpha} \int f^\alpha\,d\mu.
\en
Then $Qf \leq 0$ for the functional (6.1) in the whole class of probability densities
$f$ on $E$ having finite $L^\alpha(\mu)$-norm. According to Lemma 6.3 and using after the change 
$d = -c$, this is equivalent to the validity of the same inequality for the special function 
$f = (g + d)_+^{\beta - 1}$ where the constant $d$ is determined in a unique way by the condition
$$
\int (g + d)_+^{\beta - 1}\,d\mu = 1.
$$
In this case,
\bee
Q f 
 & = &
\int (g + d)_+^{\beta - 1}\,g\,d\mu - \frac{1}{\alpha} \int (g + d)_+^{\alpha(\beta - 1)}\,d\mu \\
 & = &
\int (g + d)_+^{\beta - 1}\,((g+d) - d)\,d\mu - \frac{1}{\alpha} \int (g + d)_+^\beta\,d\mu \\
 & = &
\Big(1 - \frac{1}{\alpha}\Big) \int (g + d)_+^\beta\,d\mu - d.
\ene
Thus, the property $Qf \leq 0$ is equivalent to $ \int (g + d)_+^\beta\,d\mu \leq \beta d$.
Since the latter integral is non-negative, necessarily $d \geq 0$, and the last relation is simplified to
\be
\int (g + d)^\beta\,d\mu \leq \beta d.
\en
Applying Jensen's inequality, we therefore obtain that
$$
\beta d \geq \Big(\int g\,d\mu + d\Big)^\beta \geq d^\beta,
$$
which can be solved as $d \leq \beta^{\frac{1}{\beta - 1}}$. Using this in (6.6), we get
\be
\int g^\beta\,d\mu \leq \beta^\alpha.
\en

In the general case where we do not assume in advance that  $g \in L^\beta(\mu)$
and start with the hypothesis (6.5), one may apply the previous step to the functions
$g_N = \min(g,N)$, $N = 1,2,\dots$. These functions are bounded, and the relation (6.5)
holds for all $g_N$. Hence we obtain the upper bound (6.7) with $g_N$ in place of $g$,
and letting $N \rightarrow \infty$ leads to (6.7) for the function $g$ itself.

Finally, starting from the inequality (6.1) with a general constant $K>0$, we may rewrite it as (6.3) 
for the function $g/(\alpha K)$ in place of $g$. Hence, by (6.7),
$$
\int g^\beta\,d\mu \leq \beta^\alpha\,(\alpha K)^\beta,
$$
which yields
$$
K \geq \frac{1}{\alpha \beta^{\alpha/\beta}}\,\|g\|_\beta.
$$
Since $ \beta^{\alpha/\beta} = \beta^{\frac{1}{\beta - 1}} < e$, the lower bound in (6.2) follows.
\qed

\vskip7mm
\section{{\bf Proof of Theorems 4.1-4.2}}
\setcounter{equation}{0}

\vskip2mm
\noindent
The relation
\be
\int g\,d\nu \, = \, \int f g\,d\mu \, \leq \, T_\alpha(\nu||\mu)
\en
for a probability measure $\nu$ with density $f = d\nu/d\mu$
may be rewritten as
$$
\int f\big(1 + (\alpha-1) g\big)\,d\mu \leq \int f^\alpha\,d\mu.
$$
It is of the form (6.1) with constant $K=1$ for the function $\widetilde g = 1 + (\alpha-1) g$ 
in place of $g$. Applying Corollary 6.2, we conclude that (7.1) may hold for any 
probability density $f$ on $(E,\mu)$ such that the integral $\int fg\,d\mu$ exists, 
only if $\tilde g_+ \in L^\beta(\mu)$. But the latter is equivalent to saying that
$g_+ \in L^\beta(\mu)$, which is thus a necesary condition for (7.1).

\vskip5mm
{\bf Remark.} In fact, Corollary 6.2 gives an additional information via the relation (6.2).
Applying it with $K=1$ to $\tilde g$, we also obtain that the necessary condition for (7.1) is
the upper moment bound
\be
\int \big(1 + (\alpha-1)g\big)_+^\beta\,d\mu \leq (e\alpha)^\beta.
\en 

\vskip3mm
{\bf Proof of Theorem 4.1.}
With the assumption that $g_+ \in L^\beta(\mu)$, it is sufficient to combine Lemmas 4.3 with 5.1
and to look at the value of the functional $R$ on the extreme density 
$$
f = \beta^{1-\beta}\, (g-c)_+^{\beta - 1},
$$ 
where $c$ is a unique solution to the equation
\be
\int (g-c)_+^{\beta - 1}\,d\mu \, = \, \beta^{\beta - 1},
\en
as indicated in (4.3). This equation corresponds to the requirement $\int f\,d\mu = 1$. 

Applying (7.3), we have
\bee
\int fg\,d\mu
 & = &
\beta^{1-\beta} \int (g-c)_+^{\beta - 1}\, \big((g-c) + c\big)\,d\mu \\
 & = &
\beta^{1-\beta} \int (g-c)_+^\beta\,d\mu \, + \, c.
\ene
Secondly, since
$$
f^\alpha = \beta^{\alpha(1-\beta)}\, (g-c)_+^{\alpha(\beta - 1)} = 
\beta^{-\beta}\, (g-c)_+^\beta,
$$
we have
\bee
\frac{1}{\alpha - 1}\,\left[\int f^\alpha\,d\mu - 1\right] 
 & = &
\frac{\beta^{-\beta}}{\alpha - 1}\,\int (g-c)_+^\beta\,d\mu \, - \,
\frac{1}{\alpha - 1} \\
 & = &
\beta^{-\beta}\,(\beta - 1)\,\int (g-c)_+^\beta\,d\mu \, - \, (\beta - 1).
\ene
Hence
$$
Rf \, = \, \beta^{-\beta} \int (g-c)_+^\beta\,d\mu \, + \, c + \beta - 1.
$$
Using this extremeizer, one may rewrite the property 
"$Rf \leq 0$ for all $f$", that is, (4.1), as
\be
\int (g-c)_+^\beta\,d\mu \ \leq \ -\beta^\beta\,(c + \beta - 1),
\en
which is the condition (4.2). This proves Theorem 4.1.
\qed

\vskip5mm
{\bf Proof of Theorem 4.2.} Suppose that (7.4) is satisfied
for some point $c_1$ in place of $c$,
\be
\int (g-c_1)_+^\beta\,d\mu \ \leq \ -\beta^\beta\,(c_1 + \beta - 1).
\en
Consider the function
$$
\psi(t) = \int (g-t)_+^\beta\,d\mu \, + \, \beta^\beta\,(t + \beta - 1), \quad t \in \R,
$$
so that (7.5) is equivalent to $\psi(c_1) \leq 0$. The function $\psi$ is strictly convex and 
is differentiable on $\R$, with $\psi(-\infty) = \psi(\infty) = \infty$. It attains minimum 
at a unique point $t$, namely -- at which
$$
\psi'(t) = \beta\, \Big[-\int (g-t)_+^{\beta - 1}\,d\mu \, + \, \beta^{\beta - 1}\Big] = 0.
$$
But this is exactly the equation (7.3), so $t = c$ for the constant from Theorem 4.1. 
Thus, if $\psi(c_1) \leq 0$, then $\psi(c) \leq 0$, so that the condition (4.2) is satisfied.
\qed

\vskip7mm
\section{{\bf Necessary and Sufficient Conditions. Moment Bounds}}
\setcounter{equation}{0}

\vskip2mm
\noindent
Since the description given in Theorem 4.1 for the property
\be
\int g\,d\nu \ \leq \ T_\alpha(\nu|\mu) \quad {\rm for \ any} \ \nu <\!\!< \mu
\en
is somewhat implicit, it would be interesting to get more tractable conditions,
necessary and sufficient, even if not simultaneously. Here we mention some of
such conditions, together with lower and upper bounds on the constant $c$
appearing in (4.2)-(4.3). To avoid situations when $D_\alpha(\nu|\mu)$ is finite,
but the integral in (8.1) does not exist, we assume that $g_+ \in L^\beta(\mu)$.

\vskip5mm
{\bf Proposition 8.1.} {\sl For $(8.1)$ to hold, it is necessary that
$\int g\,d\mu \leq 0$ and
\be
\int \Big(1 + \frac{1}{\beta-1}\, g\Big)_+^{\beta - 1}\,d\mu \leq 1,
\en
and it is sufficient that
\be
\int \Big(1 + \frac{1}{\beta}\, g\Big)_+^\beta\,d\mu \leq 1.
\en
}

\vskip3mm
{\bf Proof.} Applying (8.1) to the measure $\nu = \mu$, we get $\int g\,d\mu \leq 0$. 

The condition (8.3) is provided by Theorem 4.2, since it is exactly the inequality (4.2) 
with $c = - \beta$. Let us however describe a simple direct argument showing the sufficiency
of (8.3). By H\"older's inequality,
if the density $f$ of $\nu$ with respect to $\mu$ belongs to $L^\alpha(\mu)$, then
\bee
\int fg\,d\mu
 & = &
\beta \int f \Big(1 + \frac{1}{\beta}\, g\Big)\,d\mu - \beta \\
 & \leq &
\beta\, \Big(\int f^\alpha\,d\mu\Big)^{1/\alpha}\,
\Big(\int \Big(1 + \frac{1}{\beta}\, g\Big)_+^\beta\,d\mu\Big)^{1/\beta}  - \beta \\
 & \leq &
\beta\, \Big(\int f^\alpha\,d\mu\Big)^{1/\alpha} - \beta \\
 & \leq &
\frac{1}{\alpha - 1}\, \int f^\alpha\,d\mu - \frac{1}{\alpha - 1},
\ene
where we made use of the elementary inequality $\alpha t \leq t^\beta + \alpha - 1$
($\alpha,t \geq 1$) in the last step. Thus, we obtain (8.1).

To derive (8.2), one may just apply (8.1) to the measure $\nu$ with density
$$
f = \frac{1}{A}\,\big(1 + (\alpha-1)\,g\big)_+^{\beta - 1}, \quad
A = \int \big(1 + (\alpha-1)g\big)_+^{\beta - 1}\,d\mu.
$$
\qed

\vskip2mm
{\bf Remark.} As $\alpha \downarrow 1$, that is, $\beta \uparrow \infty$,
(8.2) and (8.3) are asymptotically optimal. In the limit they yield $\int e^g\,d\mu \leq 1$ 
which is necessary and sufficient for the relation $\int g\,d\nu \leq D(\nu|\mu)$.
Nevertheless, being quite explicit and working, (8.2)-(8.3) are not sharp enough.

\vskip2mm
Let us now consider one-sided moment bounds for $g$ under the measure $\mu$
in presence of the hypothesis such as (8.1). For this aim, we return to Theorem 4.2 and 
recall the description of (8.1) in terms of the property
\be
\int (g-c)_+^\beta\,d\mu \, \leq \, -\beta^\beta\, (c + \beta - 1),
\en
where $c$ solves the equation
\be
\int (g-c)_+^{\beta - 1}\,d\mu \, = \, \beta^{\beta- 1}.
\en

\vskip2mm
{\bf Proposition 8.2.} {\sl Under $(8.4)-(8.5)$, necessarily $c \leq - \beta$.
Furthermore, if $g \in L^1(\mu)$, then
\bee
c \geq - \beta + \int g\,d\mu & \ \ in \ the \ case \ 1 < \alpha \leq 2, \\
c \geq - 4 + \alpha \int g\,d\mu & in \ the \ case \ \alpha \geq 2.
\ene
In particular, in the corresponding cases,
\be
\int g_+^\beta\,d\mu \, \leq \, \beta^\beta\, \Big(1 - \int g\,d\mu\Big),
\en
\be
\int g_+^\beta\,d\mu \, \leq \, \beta^\beta\, \Big(4 - \alpha \int g\,d\mu\Big).
\en
}

\vskip2mm
{\bf Proof.} The weaker upper bound $c \leq -(\beta - 1)$ immediately follows 
from (8.4). To refine it, we use (8.5) and apply Jensen's inequality to get
\bee
\beta^\beta
 & = &
\left[\int (g-c)_+^{\beta - 1}\,d\mu\right]^\alpha \\
 & \leq &
\int (g-c)_+^{\alpha (\beta - 1)}\,d\mu \, = \, \int (g-c)_+^\beta\,d\mu \, \leq \, 
-\beta^\beta\,(c + \beta - 1). 
\ene
Hence, $1 \leq -(c + \beta - 1)$ proving the first statement.

For the lower bound on $c$ in the case $1 < \alpha \leq 2$, we use the convexity of
$t \rightarrow (t-c)_+^{\beta- 1}$ and apply Jensen's inequality in (8.5) 
to get
$$
\beta^{\beta- 1} \, = \, \int (g-c)_+^{\beta- 1}\,d\mu \, \geq \,
\Big(\int g\,d\mu - c\Big)_+^{\beta - 1}.
$$
Hence
$$
\beta \, \geq \, \Big(\int g\,d\mu - c\Big)_+ \, \geq \, \int g\,d\mu - c
$$
which is the first lower bound. Using this lower bound together with (8.4), we conclude that
$$
\int g_+^\beta\,d\mu \, \leq \, \int (g-c)_+^\beta\,d\mu \, \leq \, 
\beta^\beta\,(-c - \beta + 1)
 \, \leq \, \beta^\beta\, \Big(1 - \int g\,d\mu\Big)
$$
which is (8.6).

In the case $\alpha \geq 2$, from (8.4) it follows that
$\int (g-c)_+^\beta\,d\mu \leq -\beta^\beta c$. By Jensen's inequality,
$$
\int (g-c)_+^\beta\,d\mu \, \geq \, \Big(\int g\,d\mu - c\Big)_+^\beta
$$
giving
$$
\int g\,d\mu - c \, \leq \, \Big(\int g\,d\mu - c\Big)_+ \, \leq \, 
\beta (-c)^{1/\beta}.
$$
Equivalently, substituting $t = -c$, $p = \beta$, $q = \alpha$, $a = - \int g\,d\mu$, 
we arrive at the relation
$$
\varphi(t) \, \equiv \, t - p t^{1/p} \, \leq \, a.
$$
This function is convex in $t \geq 0$ and positive for $t > t_0 = p^q$, with
$\varphi(t_0) = 0$, $\varphi'(t_0) = 1 - t_0^{-1/q} = 1/q$. Hence, for all $t \geq t_0$,
$$
\varphi(t) \geq \varphi(t_0) + \varphi'(t_0) (t - t_0) = \frac{1}{q}\,(t - t_0).
$$
Once $\varphi(t) \leq a$ and $t \geq t_0$, we then get $t \leq t_0 + qa$. But 
$t_0 \leq 4$ whenever $q \geq 2$. Indeed, for the function $\psi(q) = \log t_0 = q\log p$ 
we have $\psi''(q) = \frac{1}{q(q-1)^2} > 0$, so it is convex. In addition,
$\psi'(\infty) = 0$, so it is decreasing. Hence, $\psi(q) \leq \psi(2)$ for all $q \geq 2$,
i.e., $p^q \leq 4$. This gives the required upper bound on $c$.
Again, using it in (8.4), we conclude that
\bee
\int g_+^\beta\,d\mu 
 & \leq &
\beta^\beta\, (-c - \beta + 1) \\
 & \leq &
\beta^\beta\, \Big(5 - \beta - \alpha \int g\,d\mu\Big) \ \leq \ 
\beta^\beta\, \Big(4 - \alpha \int g\,d\mu\Big)
\ene
which is (8.7).
\qed

\vskip7mm
\section{{\bf Proof of Theorem 2.3 (the case $\alpha \geq 2$)}}
\setcounter{equation}{0}

\vskip2mm
\noindent
Towards Theorem 2.3, we now turn to the question of whether or not it is possible to bound
$$
\Big|\int u\,d\nu - \int u\,d\mu\Big|
$$
in terms of the divergence $T_\alpha(\nu||\mu)$ with a fixed parameter $\alpha > 1$
uniformly over all $\nu < \!\! < \mu$ for a given function $u$.
For this aim, let us note that when $\nu = \nu_\ep$ has density of the form $f_\ep = 1 + \ep h$
with a bounded function $h$ and small $\ep$, the quantity $\Delta$ is of the
order $\ep$, while  $T_\alpha(\nu_\ep||\mu)$  is of the order $\ep^2$. Hence, it is natural
to consider the square root $\sqrt{T_\alpha(\nu||\mu)}$ as a potential upper bound on $\Delta$.

Thus, let us consider the family of inequalities of the form
\be
\Big|\int u\,d\nu\Big| \leq K\sqrt{T_\alpha(\nu||\mu)},
\en
which are required to hold for a fixed $\mu$-integrable function $u$ on $E$ with $\mu$-mean 
zero and a constant $K \geq 0$ in the class of all probability measures $\nu$ on $E$
such that the integral in (9.1) exists. Note that the assumption that $\int u\,d\mu = 0$ is necessary, 
as it follows from (9.1) in the case $\nu = \mu$. First we observe that (9.1) causes the function
$u$ to belong to the space $L^2(\mu)$. 

\vskip5mm
{\bf Lemma 9.1.} {\sl Under $(9.1)$ with $\alpha>1$, we necessarily have
\be
\int u^2\,d\mu \leq \frac{K^2 \alpha}{2}.
\en
}

\vskip3mm
{\bf Proof.}
Given a bounded measurable function $h$ on $E$ such that $\int h\,d\mu = 0$ and $\ep>0$ 
small enough, the function $f_\ep = 1 + \ep h$ represents the density of a probability 
measure $\nu = \nu_\ep$ with respect to $\mu$. In this case, (9.1) becomes
\be
\ep\, \Big|\int u h\,d\mu\Big| \leq K\sqrt{T_\alpha(\nu_\ep|\mu)}.
\en
Furthermore, by Taylor's expansion over $\ep$,
$$
T_\alpha(\nu_\ep||\mu) = \frac{1}{\alpha - 1}\,\left[\int f_\ep^\alpha\,d\mu - 1\right]
 = \frac{\alpha \ep^2}{2}\, \int h^2\,d\mu + O(\ep^3).
$$
Inserting this in (9.3) and letting $\ep \rightarrow 0$, we arrive at
$$
\Big|\int u h\,d\mu\Big| \leq K\sqrt{\frac{\alpha}{2}}\, \|h\|_2
$$
holding for any bounded $h$. But this is equivalent to 
$\|u\|_2 \leq K\sqrt{\frac{\alpha}{2}}$.
\qed

\vskip5mm
As a consequence, we arrive at the claim $b)$ in Theorem 2.3.

\vskip5mm
{\bf Corollary 9.2.} {\sl If $\alpha \geq 2$, the best constant $K$ in $(9.1)$ satisfies
\be
\sqrt{\frac{2}{\alpha}}\,\|u\|_2 \leq K \leq \|u\|_2.
\en
}

\vskip3mm
The lower bound  in (9.4) corresponds to (9.2). As for the upper bound,
the case $\alpha=2$ follows by applying Cauchy's inequality: If $f$ is density of $\nu$ with respect 
$\mu$, then
$$
\Big|\int u\,d\nu\Big| = \Big|\int u\,(f-1)\,d\mu\Big| \leq \|u\|_2\,\|f-1\|_2 = \|u\|_2\,\sqrt{T_2(\nu||\mu)},
$$
where we used the assumption that $u$ has $\mu$-mean zero.
In the general case $\alpha \geq 2$, we just use the monotonicity of $T_\alpha$ with respect to $\alpha$,
implying $T_2 \leq T_\alpha$. This proves the upper bound in (9.4).

\vskip7mm
\section{{\bf Proof of Theorem 2.3 (the Case $1 \leq \alpha \leq 2$)}}
\setcounter{equation}{0}

\vskip2mm
\noindent
Let us recall the notation
$$
K_p^{2p-2} \, = \, \sup_{r > 0} \, \Big[r^{p - 2} \int_{|u| \geq r} |u|^p\,d\mu\Big],
$$
where $u$ is a given integrable function $u$ on $E$ with $\mu$-mean zero. We need to show
that this quantity is related to the best constant $K$ in the inequality
\be
\Big|\int u\,d\nu\Big| \leq K\sqrt{T_\alpha}, \quad {\rm where} \ \ T_\alpha = T_\alpha(\nu||\mu), \ \nu <\!\! < \mu,
\en
via the two-sided bound
\be
c_\beta K_{\beta} \leq K \leq C_\beta K_{\beta}
\en
up to some constants depending on $\beta$ only.
Recall that $\beta = \frac{\alpha}{\alpha - 1}$ denotes the conjugate power.
As we will see, (10.3) hods true in the case $1 < \alpha \leq 2$, that is, for $\beta \geq 2$ with
$$
c_\beta = \frac{1}{4}\,\beta^{-\beta} \quad {\rm and} \quad C_\beta = 2\beta^\beta.
$$

First, let us explain how one can connect (10.1) to the relations of the form
\be
\int g\,d\nu \leq T_\alpha(\nu||\mu)
\en
which we discussed before.
Using the identity $\sup_\lambda\, (\lambda a - \lambda^2) = \frac{a^2}{4}$, we have
$$
\Big(\int \frac{u}{K}\,d\nu\Big)^2 \, = \,
4\, \sup_\lambda \, \Big[\lambda \int \frac{u}{K}\,d\nu - \lambda^2\Big] \, = \,
\sup_\lambda \int \Big(\frac{4}{K} \lambda u - 4\lambda^2\Big)\,d\nu.
$$
Hence, (10.1) is reduced to the inequality of the form (10.3). That is, we have:

\vskip5mm
{\bf Lemma 10.1.} {\sl Given a positive constant $K$, the relation $(10.1)$ is equivalent to
\be
\int \Big(\frac{4}{K} \lambda u - 4\lambda^2\Big)\,d\nu \, \leq \, T_\alpha(\nu|\mu)
\en
with arbitrary $\lambda \in \R$.
}

\vskip5mm
{\bf Proof of Theorem 2.3} $a)$. Thus let $1 < \alpha \leq 2$.

\vskip2mm 
{\bf Lower bound on $K$}.
First assume that $u \in L^\beta(\mu)$ has $\mu$-mean zero. Applying Proposition 8.2
 to the function $g = \frac{4}{K}\,\lambda u - 4\lambda^2$ with 
$\lambda>0$, we get the moment bound (8.6) as a necessary condition for (8.4), namely
$$
\int \Big(\frac{4}{K}\,\lambda u - 4\lambda^2\Big)_+^\beta\,d\mu \, \leq \, 
\beta^\beta (1 + 4\lambda^2).
$$
We restrict the integral to the set of points where $u \geq 2K\lambda$, on which
$\frac{4}{K} \lambda u - 4\lambda^2 \geq \frac{2}{K} \lambda u$, and then the above 
inequality yields
$$
\frac{(2\lambda)^\beta}{K^\beta}
\int_{u \geq 2K\lambda} u^\beta\,d\mu \, \leq \, \beta^\beta (1 + 4\lambda^2).
$$
To simplify, assume that $\lambda \geq 1/2$, so that $1 + 4\lambda^2 \leq8 \lambda^2$, and
then we get
$$
\lambda^{\beta - 2}
\int_{u \geq 2K\lambda} u^\beta\,d\mu \, \leq \, 8\,(K\beta/2)^\beta.
$$
Substituting $\lambda = r/(2K)$ and applying the same inequality to $-u$, we arrive at
\be
r^{\beta - 2}
\int_{|u| \geq r} |u|^\beta\,d\mu \, \leq \, 2 \beta^\beta K^{2\beta - 2},
\quad r \geq K.
\en
In the case $0 \leq r \leq K$, there is a similar obvious bound
$$
r^{\beta - 2}
\int_{r \leq |u| \leq K} |u|^\beta\,d\mu \, \leq \, K^{2\beta - 2}.
$$
Here the right-hand side is smaller than the one in (10.5). Hence, adding the
two inequalities, we arrive at
\be
r^{\beta - 2}
\int_{|u| \geq r} |u|^\beta\,d\mu \, \leq \, 4 \beta^\beta K^{2\beta - 2},
\qquad r \geq 0.
\en

On this step, the assumption $u \in L^\beta(\mu)$ can be removed: 
(10.6) can be applied to centered truncated functions 
$u_n = v_n - \int v_n\,d\mu$, where $v_n = u$ in the case $|u| \leq n$, and $v_n = \pm n$
with the sign depending on whether $u > n$ or $u < -n$. By the integrability of $u$, we have
$\int v_n\,d\mu \rightarrow \int u\,d\mu = 0$, so that
$u_n \rightarrow u$ posintwise on $E$.
Letting $n \rightarrow \infty$ and applying Fatou's lemma, we arrive at (10.6) for $u$.
This yields the left inequality in (10.4) with  $c_\beta^{-1} = 4\beta^\beta$.
Applying Lemma 10.1, we conclude that the constant $K$ in the relation (10.1) 
satisfies the lower bound in (10.2) with the factor $c_\beta$.

\vskip4mm
{\bf Upper bound on $K$.} 
By Lemma 10.1 and Theorem 4.2 applied with the same function 
$g = \frac{4}{K}\,\lambda u - 4\lambda^2$, the relation 
(10.1) holds for all $\nu <\!\!< \mu$, if and only if, for some $c \in \R$,
\be
\int \Big(\frac{4}{K}\,\lambda u - 4\lambda^2 - c\Big)_+^\beta\,d\mu \, \leq \, 
\beta^\beta\, (1 - \beta - c).
\en

{\bf Case} $0 \leq |\lambda| \leq \sqrt{\beta}$. We derive this inequality
with $c = -\beta - 4\lambda^2$, when it becomes
\be
\varphi(\lambda) \equiv \int \big(1 + \ep \lambda u\big)_+^\beta\,d\mu
 \, \leq \, 1 + 4\lambda^2, \quad \ep = \frac{4}{K\beta}.
\en
To obtain it with some $K$ independent of $\lambda$ and $u$, note that
the definition of $L = K_\beta$ implies
$$
\int |u|^\beta\,d\mu \, \leq \, 
r^\beta + \int_{|u| \geq r} |u|^\beta\,d\mu \, \leq \, r^\beta + 
\frac{L^{2\beta - 2}}{r^{\beta - 2}} \quad (r>0).
$$
Choosing here $r = L$ leads to the upper bounds 
\be
\int |u|^\beta\,d\mu \leq 2 L^\beta, \quad \int u^2\,d\mu \leq 2 L^2,
\en
where the second bound follows from the first one in view of the assumption
$\beta \geq 2$.

The function $\varphi(\lambda)$ is convex in $\lambda$ and is twice continuously 
differentiable. In addition,
$$
\varphi'(\lambda) = \beta \ep \int u \big(1 + \ep \lambda u\big)_+^{\beta - 1}\,d\mu,
$$
so that $\varphi(0) = 1$ and $\varphi'(0) = 0$. Hence, (10.8) would follow from the 
bound $\varphi''(\lambda) \leq 2$ in the region $|\lambda| \leq \sqrt{\beta}$. We have
\bee
\varphi''(\lambda) 
 & = &
\beta (\beta - 1) \int (\ep u)^2 \big(1 + \ep \lambda u\big)_+^{\beta - 2}\,d\mu \\
 & \leq &
\beta (\beta - 1) \int (\ep u)^2 \big(1 + \ep\, |\lambda u|\big)^{\beta - 2}\,d\mu.
\ene
Using $(a+b)^p \leq 2^p\, (a^p + b^p)$ ($a,b,p \geq 0$) and the assumption on the
range of $\lambda$, we get a pointwise bound
\bee
\big(1 + \ep\, |\lambda u|\big)^{\beta - 2}
 & \leq &
2^{\beta - 2}\,\big(1 + \ep^{\beta - 2} \lambda^{\beta-2}\, |u|^{\beta - 2}\big) \\
 & \leq &
\big(2\sqrt{\beta}\big)^{\beta-2}\, \big(1 + \ep^{\beta - 2}\, |u|^{\beta - 2}\big).
\ene
It then follows from (10.9) that
\bee
\varphi''(\lambda) 
 & \leq & 
\beta (\beta - 1) \, \big(2 \sqrt{\beta}\big)^{\beta-2}
\int \left((\ep u)^2 + (\ep |u|)^\beta \right)\,d\mu \\
 & \leq & 
2^{\beta-1} \beta^{\beta/2 + 1}\, \Big((\ep L)^2 + (\ep L)^\beta\Big)
 \ \leq \
2^\beta \beta^{\beta/2 + 1}\,(\ep L)^2,
\ene
where the last inequality holds true, if $\ep L \leq 1$. Choosing $K = CL$ with 
a constant $C \geq 2$, we will do have $\ep L = 4/(C\beta) \leq 1$, and then
$$
\varphi''(\lambda) \leq 2^{\beta + 4}  \beta^{\beta/2 - 1}\,\frac{1}{C^2}.
$$
Therefore, $\varphi''(\lambda) \leq 2$ as long as $C \geq 2^{(\beta + 3)/2}  \beta^{(\beta - 2)/4}$.
To simplify, here the right-hand side may be bounded from above by $2\beta^\beta$.
Hence (10.7)--(10.8) are fulfilled with $K = 2 \beta^\beta L$.

\vskip3mm
{\bf Case} $|\lambda| \geq \sqrt{\beta}$. Choosing the value
$c = -\beta - 2\lambda^2$ in (10.7), we need to show that
\be
\psi(\lambda) \equiv 
\int \Big(1 + \ep \lambda u - \frac{2\lambda^2}{\beta}\Big)_+^\beta\,d\mu
 \, \leq \, 1 + 2\lambda^2, \quad \ep = \frac{4}{K\beta}.
\en
By the assumption, $1 - \frac{2\lambda^2}{\beta} \leq - \frac{\lambda^2}{\beta}$, so,
$$
\psi(\lambda) \leq
\int \Big(\ep \lambda u - \frac{\lambda^2}{\beta}\Big)_+^\beta\,d\mu.
$$
If $|u| < |\lambda|/(\ep \beta)$, the expression $\ep \lambda u - \frac{\lambda^2}{\beta}$
is negative and therefore does not contribute to the above integral. Hence, using the definition of $L$, 
\bee
\psi(\lambda) 
 & \leq &
\int_{|u| \geq \frac{|\lambda|}{\ep \beta}} 
\Big(\ep \lambda u - \frac{\lambda^2}{\beta}\Big)^\beta\,d\mu \\
 & \leq & 
\int_{|u| \geq \frac{|\lambda|}{\ep \beta}} (\ep\, |\lambda u|)^\beta\,d\mu \, \leq \,
(\ep |\lambda|)^{\beta} \cdot \frac{L^{2\beta - 2}}{(\frac{|\lambda|}{\ep \beta})^{\beta - 2}} \\ 
 & = & 
\lambda^2 \beta^{\beta-2}\, (\ep L)^{2\beta - 2} \ = \ 
\lambda^2 \beta^{-\beta}\, \Big(\frac{4L}{K}\Big)^{2\beta-2} \ \leq \ \lambda^2,
\ene
where the last inequality holds true for the choice $K = 4L$. Thus, (10.10) is fulfilled.
\qed

\vskip7mm
\section{{\bf Pearson-Vajda distances. Proof of Theorems 2.1--2.2}}
\setcounter{equation}{0}

\vskip2mm
\noindent
Removing the condition that the function $u$ has $\mu$-mean zero,
let us give a more flexible variant of Corollary 2.4 when (2.5) is replaced with
\be
\Big|\int u\,d\nu - \int u\,d\mu \Big| \leq K\sqrt{T_\alpha}, \quad T_\alpha = T_\alpha(\nu||\mu).
\en
This relation is required to hold for a fixed $\mu$-integrable function $u$ on $E$ 
with a constant $K \geq 0$ in the class of all probability measures $\nu$ on $E$
such that the integral over $\nu$ exists.

To derive (11.1), one may apply Corollary 2.4 to the function $\widetilde u = u - c$ 
with $c = \int u \,d\mu$. Since $|c| \leq \|u\|_{2\beta^*  -2}$, we arrive at:

\vskip5mm
{\bf Corollary 11.1.} {\sl The inequality $(11.1)$ holds true with
\be
K \leq C_\beta\,\|u\|_{2\beta^* - 2}, \quad \beta^* = \max(\beta,2).
\en
Here one may take $C_\beta = 2$ for $\alpha \geq 2$ and $C_\beta = 4\,\beta^\beta$
for $\alpha<2$.
}

\vskip5mm
Of course, (11.2) can be sharpened as $K \leq C_\beta\,\|u - c\|_{2\beta^* - 2}$, 
$c = \int u\,d\mu$. In particular, for $\alpha = 2$ (11.1) becomes
$$
\Big(\int u\,d\nu - \int u\,d\mu \Big)^2 \leq \Var_\mu(u)\,T_2.
$$
In fact, this inequality may be derived by direct arguments without appealing to
Theorem 2.3. Introducing the density $f = d\nu/d\mu$, we have, by Cauchy's inequality,
\bee
\Big(\int u\,d\nu - \int u\,d\mu \Big)^2
 & = &
\Big(\int u\,(f-1)\,d\mu \Big)^2 \\
 & \leq &
\int u^2\,d\mu \int (f-1)^2\,d\mu
 \\
 & \leq &
\int u^2\,d\mu \int (f^2-1)\,d\mu \, = \, \|u\|_2^2\ T_2.
\ene
In the resulting inequality one may replace $u$ with $u-c$, which leads to the desired
statement.

A similar argument can be used for any $\alpha>1$. By H\"older's inequality,
$$
\Big|\int u\,d\nu - \int u\,d\mu \Big| \, = \, \Big|\int u\,(f-1)\,d\mu \Big| \, \leq \,
\|u\|_\beta\, \Big(\int |f-1|^\alpha\,d\mu\Big)^{1/\alpha}.
$$
Here on the right-hand side we deal with the Pearson-Vajda distance
$$
\chi_\alpha(\nu,\mu) = \int |f-1|^\alpha\,d\mu = \|f-1\|_\alpha^\alpha
$$
which coincides with the classical Pearson  $\chi^2$-distance in the case $\alpha=2$. Thus,
\be
\Big|\int u\,d\nu - \int u\,d\mu \Big| \, \leq \,
\|u\|_\beta\, \chi_\alpha^{1/\alpha}(\nu,\mu).
\en

Now, as was shown in \cite{B-C-G2}, Proposition 3.2, the Pearson-Vajda distance is related 
to the R\'enyi divergence power by means of the following relations:
$$
T_\alpha \leq \frac{1}{\alpha-1} \Big[\big(1 + \chi_\alpha^{1/\alpha}\big)^\alpha - 1 \Big]
$$
and conversely
\be
T_\alpha \geq \frac{3}{16}\,\min\big\{\chi_\alpha,\chi_\alpha^{2/\alpha}\big\} \ \ 
(1 < \alpha \leq 2), \qquad
T_\alpha \geq \alpha 3^{-\alpha} \chi_\alpha \ \ (\alpha \geq 2).
\en
If $\alpha \leq 2$, the first bound in (11.4) can be solved for $\chi_\alpha$ and gives
$$
\chi_\alpha \, \leq \,
\max\Big\{\frac{16}{3}\,T_\alpha,\Big(\frac{16}{3}\,T_\alpha\Big)^{\alpha/2}\Big\}
 \, \leq \,
\frac{16}{3}\,\max\big\{T_\alpha,T_\alpha^{\alpha/2}\big\}.
$$
Therefore, (11.3) yields:

\vskip5mm
{\bf Proposition 11.2.} {\sl If $1 < \alpha \leq 2$, for any function $u \in L^\beta(\mu)$,
\be
\Big|\int u\,d\nu - \int u\,d\mu \Big| \leq \frac{16}{3}\,\|u\|_\beta\,
\max\big\{\sqrt{T_\alpha}, T_\alpha^{1/\alpha}\big\}.
\en
If $\alpha \geq 2$, then
\be
\Big|\int u\,d\nu - \int u\,d\mu \Big| \leq \frac{1}{\alpha}\,3^\alpha\,\|u\|_\beta\,
T_\alpha^{1/\alpha}.
\en
}

\vskip5mm
In particular, if $T_\alpha \leq 1$, we get
\be
\Big|\int u\,d\nu - \int u\,d\mu \Big| \leq \frac{16}{3}\,\|u\|_\beta
\sqrt{T_\alpha(\nu||\mu)},
\en
which is better than (11.1)--(11.2), since $2\beta - 2 > \beta$ for $\beta>2$ and hence
$\|u\|_{2\beta-2} \geq \|u\|_\beta$. In addition, the constant $C_\beta$
tends to infinity as $\beta \rightarrow \infty$, while (11.6) contains
an absolute constant. On the other hand, if $T_\alpha > 1$, the resulting bound
$$
\Big|\int u\,d\nu - \int u\,d\mu \Big| \leq \frac{16}{3}\,\|u\|_\beta\,T_\alpha^{1/\alpha}
$$
is worse in comparison with (11.1)-(11.2) when the Tsallis distance is large.

\vskip5mm
{\bf Proof of Theorems 2.1--2.2.} If
$\nu$ and $\mu$ are absolutely continuous and have densities $p$ and $q$ with respect
to some measure $\lambda$ on $(E,\mathcal E)$. Then
$$
\|w (\nu - \mu)\|_{\rm TV} = \int w\,|p - q|\,d\lambda
$$
and
$$
\int u\,d\nu - \int u\,d\mu = \int (p - q)\,d\lambda.
$$
These identities do not depend on the choice of the dominating measure $\lambda$. 
Hence
$$
\|w (\nu - \mu)\|_{\rm TV} \, = \,\sup_{|u| \leq w} \Big|\int u\,d\nu - \int u\,d\mu \Big|,
$$
where the supremum is taken over all measurable functions $u$ on $E$ such that
$|u(x)| \leq w(x)$ for all $x \in E$. To get (2.2)--(2.3), it remains to apply (11.5)--(11.6)
and use $\|u\|_\beta \leq \|w\|_\beta$. With the same argument, an application of (11.1)--(11.2)
leads to the bound (2.4).
\qed

\end{document}